\documentclass[letterpaper]{article} 
\usepackage{aaai23}  
\usepackage{times}  
\usepackage{helvet}  
\usepackage{courier}  
\usepackage[hyphens]{url}  
\usepackage{graphicx} 
\usepackage{natbib}  
\usepackage{caption} 
\usepackage{cite}
\usepackage{booktabs}
\usepackage{multirow}
\usepackage{algorithm}
\usepackage{algorithmicx}
\usepackage{algpseudocode}
\usepackage{amsmath}
\usepackage{amssymb}
\floatname{algorithm}{Algorithm}

\usepackage{newfloat}
\usepackage{listings}
\usepackage{color}
\DeclareCaptionStyle{ruled}{labelfont=normalfont,labelsep=colon,strut=off} 
\floatstyle{ruled}
\newfloat{listing}{tb}{lst}{}
\floatname{listing}{Listing}
\title{Adaptive Constraint Partition based Optimization Framework for Large-scale Integer Linear Programming(Student Abstract)} 
\author{
     Huigen Ye\textsuperscript{\rm 1,\rm 2},
     Hongyan Wang\textsuperscript{\rm 1}, 
     Hua Xu\textsuperscript{\rm 1}\thanks{Corresponding author}, 
     Chengming Wang\textsuperscript{\rm 3}, 
     Yu Jiang\textsuperscript{\rm 3} 
 }
\affiliations {
     \textsuperscript{\rm 1} State Key Laboratory of Intelligent Technology and Systems, Department of Computer Science and Technology,\\ Tsinghua University, Beijing 100084, China\\
     \textsuperscript{\rm 2} School of Computer Science and Engineering, Sun Yat-sen University, Guangzhou, China\\
     \textsuperscript{\rm 3} Meituan Inc., Block F\&G, Wangjing International R\&D Park, No.6 Wang Jing East Rd, Chaoyang District, \\ Beijing, 100102, China\\
     yehg5@mail2.sysu.edu.cn, why17@mails.tsinghua.edu.cn, xuhua@tsinghua.edu.cn,\\ wangchengming02@meituan.com, jiangyu31@meituan.com
 }

\usepackage{bibentry}

\begin{document}

\maketitle

\begin{abstract}

Integer programming problems (IPs) are challenging to be solved efficiently due to the NP-hardness, especially for large-scale IPs. To solve this type of IPs, Large neighborhood search (LNS) uses an initial feasible solution and iteratively improves it by searching a large neighborhood around the current solution. However, LNS easily steps into local optima and ignores the correlation between variables to be optimized, leading to compromised performance. This paper presents a general adaptive constraint partition-based optimization framework (ACP) for large-scale IPs that can efficiently use any existing optimization solver as a subroutine. Specifically, ACP first randomly partitions the constraints into blocks, where the number of blocks is adaptively adjusted to avoid local optima. Then, ACP uses a subroutine solver to optimize the decision variables in a randomly selected block of constraints to enhance the variable correlation. ACP is compared with LNS framework with different subroutine solvers on four IPs and a real-world IP. The experimental results demonstrate that in specified wall-clock time ACP shows better performance than SCIP and Gurobi.

\end{abstract}

\section{Introduction}
Integer programming problems (IPs) are usually NP-hard, and their algorithm design is challenging and research-worthy. Even with advances, the performance of traditional tree algorithms for IPs, such as branch-and-bound \cite{zarpellon2021parameterizing} and branch-and-cut \cite{huang2022learning}, decrease severely in high-dimensional search spaces. Therefore, Large neighborhood search (LNS) \cite{shaw1998using, pisinger2010large, NEURIPS2020_e769e03a},  defining neighborhoods and optimizing blocks that contain a subset of decision variables, has been widely used and achieved good results for many real-world large-scale IPs\cite{sonnerat2021learning,li2022mapf}. However, for IPs with millions of decision variables, the traditional LNS frameworks ignore the correlation between variables and easily step into local optima due to the fixed number of blocks. This paper proposes ACP to address the above two issues. The key idea of ACP is that it adaptively updates the number of blocks to avoid local optima, and uses a two-step variable selection method to enhance the correlation between variables to be optimized. Results on four large-scale IPs and a real-world IP verify the effectiveness of ACP.


\section{Method}
As shown in Algorithm 1, ACP starts with an initial feasible solution that is artificially constructed, and iteratively improves the current solution. Then, constraints are randomly partitioned into disjoint blocks (Step 2). Each iteration only considers one block and the variables in this block are optimized, while other variables are fixed with the value of the current optimal feasible solution (Step 3-5). The number of blocks is updated according to the objective value improvement of recent two iterations with a pre-set threshold (Step 6-7). Based on ACP, the ACP2 framework is derived that uses the subroutine solver instead of artificially constructing to generate an initial feasible solution.

\begin{algorithm}[h]
\label{alg1}
\small
\caption{The framework of ACP} 
\hspace*{0.02in} {\textbf{Input:}} 
An IP $P$ with constraints $C$ and variables $X$, an initial feasible solution $S_X$, a subroutine solver $F$  \\     
\hspace*{0.02in} {\textbf{Output:}} 
A solution $S_X$                                             
\begin{algorithmic}[1]
\While{Specified wall-clock time not met}
    \State $C = C_1 \cup C_2 \cup \dots \cup C_{k - 1} \cup C_k$ \Comment{\textit{Constraint Partition}}
	\State $X_{sub} \gets$ variables selected in a random block $C_i$
	\State $X_{opt} \gets S_X$
	\State $S_X \gets \text{FIX\_OPTIMIZE(} P, S_X, X_{sub}, F \text{)}$\Comment{\textit{Optimization}}
	\If{$f(S_X)-f(X_{opt})<\epsilon f(X_{opt})$} \Comment{\textit{Block Update}} 
		\State $k \gets k-1$
	\EndIf
\EndWhile
\State \Return $S_X$
\end{algorithmic}
\end{algorithm}

\setlength{\parindent}{0cm}
\textbf{Constraint Partition.} For an IP $P$ with the decision variable set $X$ and the constraint set $C$, ACP randomly divides $C$ into $k$ disjoint blocks $C_1, C_2, \dots, C_k$ where $C = C_1\cup C_2 \cup \dots \cup C_k$ and $C_i \cap C_j=\emptyset$, $i\neq j$, $i,j \in [1,\dots,k]$. $k$ is a hyperparameter that is different for different IPs.

\textbf{Optimization.} At each iteration, one block $C_i$($ i \in [1,\dots,k]$) is randomly selected, and $C_i$ differs at different iterations. Given a feasible solution $S_x$ of the current IP $P$, all decision variables $X_{sub}$ in the randomly selected subset $C_i$ are treated as the local neighborhood of the search. Then in function FIX\_OPTIMIZE, a subroutine solver $F$ (SCIP, Gurobi) is used to search the optimal solution of a sub-IP with decision variable set $X_{sub}$. All integer variables in $X\setminus X_{sub}$ are fixed with the value of the current optimal feasible solution $X_{opt}$. The new solution $S_X$ is obtained by recombining $X_{sub}$ and $X\setminus X_{sub}$.

\textbf{Block update.} To avoid getting stuck in local optima, ACP updates the number of blocks $k$ with an optimization threshold of objective value improvement $\epsilon$. If the improvement rate $f(S_X)-f(X_{opt})$ of the objective function $f(S_X)$ after one iteration is less than $\epsilon$, the current solution is likely to be a local optimum. In this case, ACP reduces the number of blocks $k$ to expand the neighborhood to jump out of the local optimum. Additional details such as the initial block number $k$ and optimization threshold $\epsilon$ are given in the appendix.

\section{Experiments}
\textbf{Experimental Settings.}
With different subroutine solvers, we compare ACP with LNS framework \cite{NEURIPS2020_e769e03a} and the original subroutine solver on maximum independent set (IS), minimum vertex cover (MVC), maximum cut (MAXCUT), minimum set covering(SC) and one real-world large-scale IP in the internet domain. For IS, MVC and MAXCUT, we use random graphs of 1,000,000 points and 3,000,000 edges. For SC, we use random problem of 1,000,000 items and 1,000,000 sets. For the real-world large-scale IP, it has more than 800,000 decision variables and 50,000 constraints. All experiments
are repeated 5 times and the average of the metric is recorded.

\textbf{Results and Analysis.}
Table 1 shows the results of all the related methods on five IPs, and the best results are in bold fonts. It can be seen that ACP reliably offers remarkable improvements over LNS framework with different subroutine solvers and the original subroutine solver. Compared with SCIP, ACP obviously achieves much better performance, especially for MAXCUT, SC and the real-world IP. Although the improvement compared with Gurobi is lower than that of SCIP, ACP still outperforms Gurobi on all the five IPs. We also find that with block partition and adaptive update of block number, ACP improves the performance of LNS.
\begin{figure}[t]
\centering
    \includegraphics[scale=0.1]{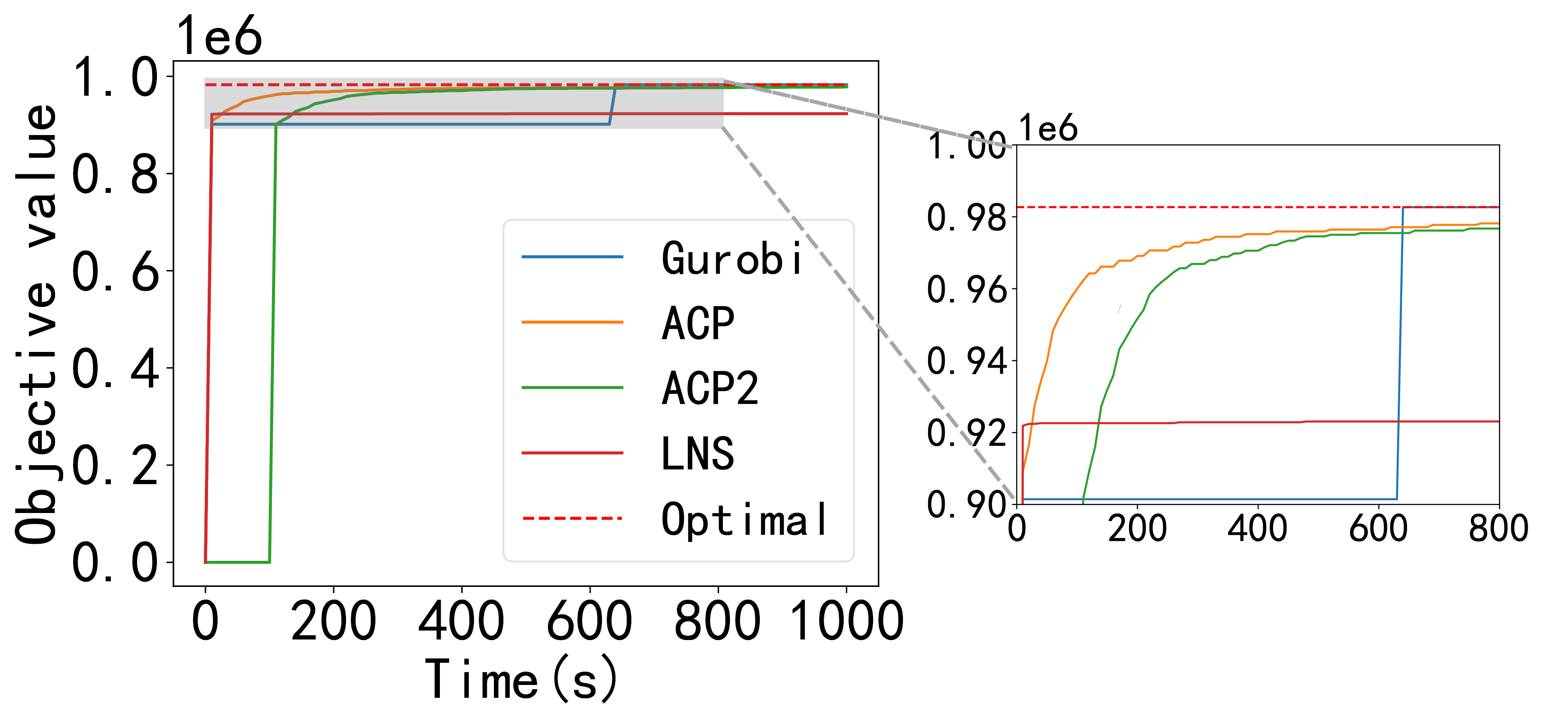}
\caption{The objective value variation for real-world IP.}
\end{figure}
For clearer comparison, we further show the objective value variation with running time for the real-world large-scale IP in Figure 1. It can be found that ACP exhibits noteworthy advantages within limited wall-clock time, which can find the approximate optimal solution faster than other methods, even for the fastest commercial solver Gurobi. Due to the length limitation, we show other comparison results of all the related methods on small, medium and large-scale IPs in the appendix. The results also demonstrate the obvious advantages of ACP over other baseline methods. %

\begin{table}[t]
\tiny
\begin{tabular}{|c|l|l|l|l|}
\hline
\textbf{Problem}                                                                    & \textbf{\begin{tabular}[c]{@{}c@{}}SCIP \\ based\end{tabular}} & \textbf{Objective value}          & \textbf{\begin{tabular}[c]{@{}c@{}}Gurobi\\ based\end{tabular}} & \textbf{Objective value}         \\ \hline
\multirow{4}{*}{\begin{tabular}[c]{@{}c@{}}IS\\ (Maximize)\end{tabular}}            & SCIP                                                           & 7866.55 $\pm$ 55.47               & Gurobi                                                          & 215365.65 $\pm$ 113.5            \\
                                                                                    & LNS                                                            & 194733.92 $\pm$ 79.17             & LNS                                                             & 220368.34 $\pm$ 557.89           \\
                                                                                    & ACP                                                            & \textbf{207901.54 $\pm$ 1531.18}  & ACP                                                             & 227559.62 $\pm$ 57.26            \\
                                                                                    & ACP2                                                           & 196742.1 $\pm$ 116.6              & ACP2                                                            & \textbf{227636.51 $\pm$ 432.18}  \\ \hline
\multirow{4}{*}{\begin{tabular}[c]{@{}c@{}}MVC\\ (Minimize)\end{tabular}}           & SCIP                                                           & 490857.44 $\pm$ 164.24            & Gurobi                                                          & 283170.59 $\pm$ 317.06           \\
                                                                                    & LNS                                                            & 304860.39 $\pm$ 313.2             & LNS                                                             & 277870.09 $\pm$ 300.8            \\
                                                                                    & ACP                                                            & \textbf{291226.0 $\pm$ 1189.94}   & ACP                                                             & 271163.92 $\pm$ 337.0            \\
                                                                                    & ACP2                                                           & 301836.42 $\pm$ 70.59             & ACP2                                                            & \textbf{271087.22 $\pm$ 295.85}  \\ \hline
\multirow{4}{*}{\begin{tabular}[c]{@{}c@{}}MAXCUT\\ (Maximize)\end{tabular}}        & SCIP                                                           & 9.02 $\pm$ 1.24                   & Gurobi                                                          & 971138.1 $\pm$ 308.18            \\
                                                                                    & LNS                                                            & 553840.66 $\pm$ 45172.75          & LNS                                                             & 910662.36 $\pm$ 516.31           \\
                                                                                    & ACP                                                            & \textbf{829597.17 $\pm$ 46331.3}  & ACP                                                             & 1050400.97 $\pm$ 106.3           \\
                                                                                    & ACP2                                                           & 747434.19 $\pm$ 48697.78          & ACP2                                                            & \textbf{1053583.50 $\pm$ 473.43} \\ \hline
\multirow{4}{*}{\begin{tabular}[c]{@{}c@{}}SC\\ (Minimize)\end{tabular}}            & SCIP                                                           & 919264.06 $\pm$ 607.75            & Gurobi                                                          & 320034.15 $\pm$ 208.69           \\
                                                                                    & LNS                                                            & 584210.54 $\pm$ 158940.18         & LNS                                                             & 224435.73 $\pm$ 1253.96          \\
                                                                                    & ACP                                                            & \textbf{392082.53 $\pm$ 10206.15} & ACP                                                             & 208062.73 $\pm$ 5582.97          \\
                                                                                    & ACP2                                                           & 432368.51 $\pm$ 615.53            & ACP2                                                            & \textbf{201034.73 $\pm$ 196.14}  \\ \hline
\multirow{4}{*}{\begin{tabular}[c]{@{}c@{}}Read-world IP\\ (Maximize)\end{tabular}} & SCIP                                                           & 0.0 $\pm$ 0.0                     & Gurobi                                                          & 903119.36 $\pm$ 1588.62          \\
                                                                                    & LNS                                                            & 904387.43 $\pm$ 1560.34           & LNS                                                             & 924886.13 $\pm$ 1670.82          \\
                                                                                    & ACP                                                            & \textbf{909113.41 $\pm$ 2030.01}  & ACP                                                             & \textbf{948461.41 $\pm$ 4385.63} \\
                                                                                    & ACP2                                                           & 907155.48 $\pm$ 2199.08           & ACP2                                                            & 942235.76 $\pm$ 1861.36          \\ \hline
\end{tabular}
\caption{Comparison results for ACP and LNS with different subroutine solvers on different benchmark IPs.}
\end{table}

\section{Conclusion}
This paper presents ACP for large-scale IPs that can efficiently use any existing optimization solver as a subroutine. The results and analysis show that the ACP framework is obviously superior to the mainstream method in specified wall-clock time. In the future, we will study the combination of LNS and graph neural networks for large-scale IPs.
\section{Acknowledgement}
This research was supported by Meituan.

{\small \bibliography{aaai23}}



\appendix
    \part*{-Appendix-}

\setcounter{table}{0}
\section{Experimental Settings}
This section shows the experimental settings used in our paper, including the description of all the five IP benchmark problems and implementation details. 
\subsection{Problem Description}
To study the performance of ACP, we use four classical integer programming problems (IPs), including maximum independent set (IS), minimum vertex cover (MVC), maximum cut (MAXCUT), minimum set covering (SC), and a real-world problem represented as large-scale IP in this paper. Each of the five IPs contains three kinds of scales, small, medium and large, which include ten thousands of, hundred thousands of and millions of decision variables, respectively. The details of all the IP benchmarks are shown in Table 1. The specified time in Table 1 are the total pre-specified running times for each method on each IP, i.e., the wall-clock time in Algorithm 1.



\begin{table}[h]
\tiny
\begin{tabular}{|c|c|c|c|}
\hline
\textbf{Scale}                   & \textbf{Problem}       & \textbf{Description}                                                                                             & \multicolumn{1}{l|}{\textbf{Specified time}} \\ \hline
\multirow{5}{*}{Small}  & IS            & 10,000 nodes and 30,000 edges                                                                             & 10s                                  \\
                        & MVC           & 10,000 nodes and 30,000 edges                                                                             & 10s                                  \\
                        & MAXCUT        & 10,000 nodes and 30,000 edges                                                                             & 10s                                  \\
                        & SC            & \begin{tabular}[c]{@{}c@{}}20,000 items and 20,000 set where\\ each item appears in 4 sets\end{tabular}   & 20s                                  \\
                        & Real-world IP & 10,000 decision variables                                                                 & 10s                                  \\ \hline
\multirow{5}{*}{Medium} & IS            & 100,000 nodes and 300,000 edges                                                                           & 100s                                  \\
                        & MVC           & 100,000 nodes and 300,000 edges                                                                           & 100s                                  \\
                        & MAXCUT        & 100,000 nodes and 300,000 edges                                                                           & 180s                                  \\
                        & SC            & \begin{tabular}[c]{@{}c@{}}200,000 items and 200,000 set where\\ each item appears in 4 sets\end{tabular} & 100s                                  \\
                        & Real-world IP & 100,000 decision variables                                                             & 50s                                  \\ \hline
\multirow{5}{*}{Large}  & IS            & 1000,000 nodes and 3000,000 edges                                                                         & 1500s                                  \\
                        & MVC           & 1000,000 nodes and 3000,000 edges                                                                         & 1500s                                  \\
                        & MAXCUT        & 1000,000 nodes and 3000,000 edges                                                                         & 1800s                                  \\
                        & SC            & \begin{tabular}[c]{@{}c@{}}200,000 items and 200,000 set where\\ each item appears in 4 sets\end{tabular} & 500s                                  \\
                        & Real-world IP & 1000,000 decision variables                                                                          & 100s                                  \\ \hline
\end{tabular}
\caption{Details of five IP benchmarks.}
\end{table}

\subsection{Implementation details}

\setlength{\parindent}{0cm}

For a fair comparison, all the related methods are run 5 times with 5 different random seeds on each IP on the same device, Intel Core i9-10900K processor with 10 cores, 2 NVIDIA GeForce RTX 3090 GPU card with 24GB GPU memory. The mean of the results is recorded. For different scales of IP benchmarks, we use different initial numbers of blocks $k$. Table 2 gives all the $k$ cases for IS, MVC, MAXCUT, SC and real-world IP. The $\epsilon$ used to determine the objective value improvement is set to different values for different scales of IPs, which is shown in Table 3. 

In the code implementation, to adjust the frequency of dynamic block number changes, we set $t$, which means that when the improvement rates of consecutive $t$ iterations are all less than the threshold $\epsilon$, ACP will reduce the number of blocks $k$ to expand the neighborhood to jump out of the local optimum, as shown in Table 4. To adapt to different problem sizes, we also set $p$, which means the proportion of the maximum running time of each iteration to the total time, as shown in Table 5. 


~\\
~\\
\begin{table}[H]
\tiny
\begin{tabular}{|c|c|cccccc|}
\hline
\multirow{3}{*}{Scale}  & \multirow{3}{*}{Problem} & \multicolumn{6}{c|}{Initial block number $k$}                                                                                                \\ \cline{3-8} 
                        &                          & \multicolumn{3}{c|}{SCIP}                                                       & \multicolumn{3}{c|}{Gurobi}                                \\ \cline{3-8} 
                        &                          & \multicolumn{1}{c|}{LNS} & \multicolumn{1}{c|}{ACP} & \multicolumn{1}{c|}{ACP2} & \multicolumn{1}{c|}{LNS} & \multicolumn{1}{c|}{ACP} & ACP2 \\ \hline
\multirow{5}{*}{Small}  & IS                       & \multicolumn{1}{c|}{2}   & \multicolumn{1}{c|}{6}   & \multicolumn{1}{c|}{4}    & \multicolumn{1}{c|}{2}   & \multicolumn{1}{c|}{4}   & 4    \\
                        & MVC                      & \multicolumn{1}{c|}{3}   & \multicolumn{1}{c|}{5}   & \multicolumn{1}{c|}{4}    & \multicolumn{1}{c|}{2}   & \multicolumn{1}{c|}{4}   & 4    \\
                        & MAXCUT                   & \multicolumn{1}{c|}{2}   & \multicolumn{1}{c|}{3}   & \multicolumn{1}{c|}{2}    & \multicolumn{1}{c|}{2}   & \multicolumn{1}{c|}{3}   & 3    \\
                        & SC                       & \multicolumn{1}{c|}{4}    & \multicolumn{1}{c|}{10}    & \multicolumn{1}{c|}{10}     & \multicolumn{1}{c|}{3}   & \multicolumn{1}{c|}{7}   & 8    \\
                        & Real-world IP            & \multicolumn{1}{c|}{5}    & \multicolumn{1}{c|}{10}    & \multicolumn{1}{c|}{10}     & \multicolumn{1}{c|}{2}    & \multicolumn{1}{c|}{3}    &  3    \\ \hline
\multirow{5}{*}{Medium} & IS                       & \multicolumn{1}{c|}{6}   & \multicolumn{1}{c|}{8}   & \multicolumn{1}{c|}{8}    & \multicolumn{1}{c|}{3}   & \multicolumn{1}{c|}{6}   & 6    \\
                        & MVC                      & \multicolumn{1}{c|}{6}   & \multicolumn{1}{c|}{8}   & \multicolumn{1}{c|}{8}    & \multicolumn{1}{c|}{3}   & \multicolumn{1}{c|}{6}   & 6    \\
                        & MAXCUT                   & \multicolumn{1}{c|}{4}   & \multicolumn{1}{c|}{6}   & \multicolumn{1}{c|}{6}    & \multicolumn{1}{c|}{4}   & \multicolumn{1}{c|}{5}   & 4    \\
                        & SC                       & \multicolumn{1}{c|}{6}    & \multicolumn{1}{c|}{12}    & \multicolumn{1}{c|}{12}     & \multicolumn{1}{c|}{3}   & \multicolumn{1}{c|}{5}   & 6    \\
                        & Real-world IP            & \multicolumn{1}{c|}{6}    & \multicolumn{1}{c|}{10}    & \multicolumn{1}{c|}{10}     & \multicolumn{1}{c|}{3}    & \multicolumn{1}{c|}{4}  & 4    \\ \hline
\multirow{5}{*}{Large}  & IS                       & \multicolumn{1}{c|}{8}   & \multicolumn{1}{c|}{10}  & \multicolumn{1}{c|}{10}   & \multicolumn{1}{c|}{6}   & \multicolumn{1}{c|}{6}   & 6    \\
                        & MVC                      & \multicolumn{1}{c|}{8}   & \multicolumn{1}{c|}{10}  & \multicolumn{1}{c|}{10}   & \multicolumn{1}{c|}{3}   & \multicolumn{1}{c|}{8}   & 8    \\
                        & MAXCUT                   & \multicolumn{1}{c|}{15}  & \multicolumn{1}{c|}{15}  & \multicolumn{1}{c|}{15}   & \multicolumn{1}{c|}{6}   & \multicolumn{1}{c|}{10}  & 5    \\
                        & SC                       & \multicolumn{1}{c|}{20}  & \multicolumn{1}{c|}{25}  & \multicolumn{1}{c|}{25}   & \multicolumn{1}{c|}{4}   & \multicolumn{1}{c|}{5}   & 5    \\
                        & Real-world IP            & \multicolumn{1}{c|}{6}   & \multicolumn{1}{c|}{50}  & \multicolumn{1}{c|}{50}   & \multicolumn{1}{c|}{3}   & \multicolumn{1}{c|}{7}   & 6    \\ \hline
\end{tabular}
\caption{The initial block number $k$ of all IPs with small, medium and large scales.}
\end{table}

~\\

\begin{center}
\begin{table}[H]
\centering
\tiny
\begin{tabular}{|c|c|cccc|}
\hline
\multirow{3}{*}{Scale}  & \multirow{3}{*}{Problem} & \multicolumn{4}{c|}{Optimization threshold $\epsilon$}                                          \\ \cline{3-6} 
                        &                          & \multicolumn{2}{c|}{SCIP}                                 & \multicolumn{2}{c|}{Gurobi}         \\ \cline{3-6} 
                        &                          & \multicolumn{1}{c|}{ACP}    & \multicolumn{1}{c|}{ACP2}   & \multicolumn{1}{c|}{ACP}   & ACP2   \\ \hline
\multirow{5}{*}{Small}  & IS                       & \multicolumn{1}{c|}{0.002}  & \multicolumn{1}{c|}{0.01}   & \multicolumn{1}{c|}{0.01}  & 0.01   \\
                        & MVC                      & \multicolumn{1}{c|}{0.002}  & \multicolumn{1}{c|}{0.01}   & \multicolumn{1}{c|}{0.01}  & 0.01   \\
                        & MAXCUT                   & \multicolumn{1}{c|}{0.1}    & \multicolumn{1}{c|}{0.01}   & \multicolumn{1}{c|}{0.01}  & 0.01   \\
                        & SC                       & \multicolumn{1}{c|}{0.002}  & \multicolumn{1}{c|}{0.01}   & \multicolumn{1}{c|}{0.01}  & 0.01   \\
                        & Real-world IP            & \multicolumn{1}{c|}{0.001}  & \multicolumn{1}{c|}{0.05}   & \multicolumn{1}{c|}{0.001} & 0.001 \\ \hline
\multirow{5}{*}{Medium} & IS                       & \multicolumn{1}{c|}{0.1}    & \multicolumn{1}{c|}{0.01}   & \multicolumn{1}{c|}{0.01}  & 0.01   \\
                        & MVC                      & \multicolumn{1}{c|}{0.002}  & \multicolumn{1}{c|}{0.01}   & \multicolumn{1}{c|}{0.01}  & 0.01   \\
                        & MAXCUT                   & \multicolumn{1}{c|}{0.1}    & \multicolumn{1}{c|}{0.01}   & \multicolumn{1}{c|}{0.005} & 0.005  \\
                        & SC                       & \multicolumn{1}{c|}{0.002}  & \multicolumn{1}{c|}{0.01}   & \multicolumn{1}{c|}{0.01}  & 0.01   \\
                        & Real-world IP            & \multicolumn{1}{c|}{0.02}   & \multicolumn{1}{c|}{0.05}   & \multicolumn{1}{c|}{0.01}  & 0.0007 \\ \hline
\multirow{5}{*}{Large}  & IS                       & \multicolumn{1}{c|}{0.1}    & \multicolumn{1}{c|}{0.01}   & \multicolumn{1}{c|}{0.01}  & 0.01   \\
                        & MVC                      & \multicolumn{1}{c|}{0.002}  & \multicolumn{1}{c|}{0.01}   & \multicolumn{1}{c|}{0.01}  & 0.01   \\
                        & MAXCUT                   & \multicolumn{1}{c|}{0.1}    & \multicolumn{1}{c|}{0.01}   & \multicolumn{1}{c|}{0.1}   & 0.005  \\
                        & SC                       & \multicolumn{1}{c|}{0.002}  & \multicolumn{1}{c|}{0.01}   & \multicolumn{1}{c|}{0.01}  & 0.01   \\
                        & Real-world IP            & \multicolumn{1}{c|}{0.0007} & \multicolumn{1}{c|}{0.0007} & \multicolumn{1}{c|}{0.01}  & 0.0007 \\ \hline
\end{tabular}
\caption{The optimization threshold $\epsilon$ of all IPs with small, medium and large scales.}
\end{table}
\end{center}

\begin{center}
\begin{table}[H]
\centering
\tiny
\begin{tabular}{|c|c|cccc|}
\hline
\multirow{3}{*}{Scale}  & \multirow{3}{*}{Problem} & \multicolumn{4}{c|}{Maximum iterations $t$}                                          \\ \cline{3-6} 
                        &                          & \multicolumn{2}{c|}{SCIP}                                 & \multicolumn{2}{c|}{Gurobi}         \\ \cline{3-6} 
                        &                          & \multicolumn{1}{c|}{ACP}    & \multicolumn{1}{c|}{ACP2}   & \multicolumn{1}{c|}{ACP}   & ACP2   \\ \hline
\multirow{5}{*}{Small}  & IS                       & \multicolumn{1}{c|}{3}  & \multicolumn{1}{c|}{3}   & \multicolumn{1}{c|}{3}  & 3   \\
                        & MVC                      & \multicolumn{1}{c|}{3}  & \multicolumn{1}{c|}{3}   & \multicolumn{1}{c|}{3}  & 3   \\
                        & MAXCUT                   & \multicolumn{1}{c|}{3}    & \multicolumn{1}{c|}{3}   & \multicolumn{1}{c|}{3}  & 3   \\
                        & SC                       & \multicolumn{1}{c|}{3}  & \multicolumn{1}{c|}{3}   & \multicolumn{1}{c|}{3}  & 3   \\
                        & Real-world IP            & \multicolumn{1}{c|}{3}  & \multicolumn{1}{c|}{2}   & \multicolumn{1}{c|}{3} & 3 \\ \hline
\multirow{5}{*}{Medium} & IS                       & \multicolumn{1}{c|}{3}    & \multicolumn{1}{c|}{3}   & \multicolumn{1}{c|}{3}  & 3   \\
                        & MVC                      & \multicolumn{1}{c|}{3}  & \multicolumn{1}{c|}{3}   & \multicolumn{1}{c|}{3}  & 3  \\
                        & MAXCUT                   & \multicolumn{1}{c|}{3}    & \multicolumn{1}{c|}{3}   & \multicolumn{1}{c|}{3} & 2  \\
                        & SC                       & \multicolumn{1}{c|}{3}  & \multicolumn{1}{c|}{3}   & \multicolumn{1}{c|}{3}  & 3   \\
                        & Real-world IP            & \multicolumn{1}{c|}{2}   & \multicolumn{1}{c|}{3}   & \multicolumn{1}{c|}{3}  & 3 \\ \hline
\multirow{5}{*}{Large}  & IS                       & \multicolumn{1}{c|}{3}    & \multicolumn{1}{c|}{3}   & \multicolumn{1}{c|}{3}  & 3  \\
                        & MVC                      & \multicolumn{1}{c|}{3}  & \multicolumn{1}{c|}{3}   & \multicolumn{1}{c|}{3}  & 3   \\
                        & MAXCUT                   & \multicolumn{1}{c|}{3}    & \multicolumn{1}{c|}{3}   & \multicolumn{1}{c|}{2}   & 2  \\
                        & SC                       & \multicolumn{1}{c|}{3}  & \multicolumn{1}{c|}{3}   & \multicolumn{1}{c|}{3}  & 3   \\
                        & Real-world IP            & \multicolumn{1}{c|}{5} & \multicolumn{1}{c|}{5} & \multicolumn{1}{c|}{3}  & 3 \\ \hline
\end{tabular}
\caption{The maximum consecutive iterations that are all less than the threshold.}
\end{table}
\end{center}

\begin{table}[H]
\tiny
\begin{tabular}{|c|c|cccccc|}
\hline
\multirow{3}{*}{Scale}  & \multirow{3}{*}{Problem} & \multicolumn{6}{c|}{Maximum proportion $p$}                                                                                                \\ \cline{3-8} 
                        &                          & \multicolumn{3}{c|}{SCIP}                                                       & \multicolumn{3}{c|}{Gurobi}                                \\ \cline{3-8} 
                        &                          & \multicolumn{1}{c|}{LNS} & \multicolumn{1}{c|}{ACP} & \multicolumn{1}{c|}{ACP2} & \multicolumn{1}{c|}{LNS} & \multicolumn{1}{c|}{ACP} & ACP2 \\ \hline
\multirow{5}{*}{Small}  & IS                       & \multicolumn{1}{c|}{0.1}   & \multicolumn{1}{c|}{0.1}   & \multicolumn{1}{c|}{0.2}    & \multicolumn{1}{c|}{0.1}   & \multicolumn{1}{c|}{0.1}   & 0.2    \\
                        & MVC                      & \multicolumn{1}{c|}{0.1}   & \multicolumn{1}{c|}{0.1}   & \multicolumn{1}{c|}{0.2}    & \multicolumn{1}{c|}{0.1}   & \multicolumn{1}{c|}{0.1}   & 0.2    \\
                        & MAXCUT                   & \multicolumn{1}{c|}{0.3}   & \multicolumn{1}{c|}{0.2}   & \multicolumn{1}{c|}{0.3}    & \multicolumn{1}{c|}{0.1}   & \multicolumn{1}{c|}{0.1}   & 0.2    \\
                        & SC                       & \multicolumn{1}{c|}{0.2}    & \multicolumn{1}{c|}{0.2}    & \multicolumn{1}{c|}{0.2}     & \multicolumn{1}{c|}{0.1}   & \multicolumn{1}{c|}{0.1}   & 0.2    \\
                        & Real-world IP            & \multicolumn{1}{c|}{0.2}    & \multicolumn{1}{c|}{0.25}    & \multicolumn{1}{c|}{0.3}     & \multicolumn{1}{c|}{0.2}    & \multicolumn{1}{c|}{0.3}    &  0.1    \\ \hline
\multirow{5}{*}{Medium} & IS                       & \multicolumn{1}{c|}{0.1}   & \multicolumn{1}{c|}{0.1}   & \multicolumn{1}{c|}{0.2}    & \multicolumn{1}{c|}{0.1}   & \multicolumn{1}{c|}{0.1}   & 0.2    \\
                        & MVC                      & \multicolumn{1}{c|}{0.1}   & \multicolumn{1}{c|}{0.1}   & \multicolumn{1}{c|}{0.2}    & \multicolumn{1}{c|}{0.1}   & \multicolumn{1}{c|}{0.1}   & 0.2    \\
                        & MAXCUT                   & \multicolumn{1}{c|}{0.1}   & \multicolumn{1}{c|}{0.1}   & \multicolumn{1}{c|}{0.1}    & \multicolumn{1}{c|}{0.1}   & \multicolumn{1}{c|}{0.1}   & 0.1    \\
                        & SC                       & \multicolumn{1}{c|}{0.2}    & \multicolumn{1}{c|}{0.2}    & \multicolumn{1}{c|}{0.2}     & \multicolumn{1}{c|}{0.1}   & \multicolumn{1}{c|}{0.1}   & 0.2    \\
                        & Real-world IP            & \multicolumn{1}{c|}{0.2}    & \multicolumn{1}{c|}{0.2}    & \multicolumn{1}{c|}{0.2}     & \multicolumn{1}{c|}{0.2}    & \multicolumn{1}{c|}{0.2}  & 0.1    \\ \hline
\multirow{5}{*}{Large}  & IS                       & \multicolumn{1}{c|}{0.1}   & \multicolumn{1}{c|}{0.1}  & \multicolumn{1}{c|}{0.1}   & \multicolumn{1}{c|}{0.2}   & \multicolumn{1}{c|}{0.1}   & 0.2    \\
                        & MVC                      & \multicolumn{1}{c|}{0.1}   & \multicolumn{1}{c|}{0.1}  & \multicolumn{1}{c|}{0.1}   & \multicolumn{1}{c|}{0.1}   & \multicolumn{1}{c|}{0.1}   & 0.2    \\
                        & MAXCUT                   & \multicolumn{1}{c|}{0.1}  & \multicolumn{1}{c|}{0.1}  & \multicolumn{1}{c|}{0.1}   & \multicolumn{1}{c|}{0.1}   & \multicolumn{1}{c|}{0.3}  & 0.2    \\
                        & SC                       & \multicolumn{1}{c|}{0.2}  & \multicolumn{1}{c|}{0.2}  & \multicolumn{1}{c|}{0.2}   & \multicolumn{1}{c|}{0.3}   & \multicolumn{1}{c|}{0.3}   & 0.2    \\
                        & Real-world IP            & \multicolumn{1}{c|}{0.1}   & \multicolumn{1}{c|}{0.1}  & \multicolumn{1}{c|}{0.1}   & \multicolumn{1}{c|}{0.2}   & \multicolumn{1}{c|}{0.2}   & 0.1    \\ \hline
\end{tabular}
\caption{The proportion of the maximum running time of each iteration to the total time.}
\end{table}

\section{Results on Different Scales of IPs}
This section shows additional results of all the related methods on five IP benchmarks with small-scale, medium-scale and large-scale, which are shown in Table 6-8, respectively. The best results are in marked bold fonts. 

Specifically, compared with SCIP and Gurobi, ACP framework can explore more solution spaces due to the idea of constraint partition. ACP gets better optimization results on all IPs in specified wall-clock time. Compared with the LNS framework, except for the five IPs of MAXCUT, due to the weak correlation between decision variables, the improvement of ACP compared with LNS is not as huge as compared with SCIP and Gurobi. However, in the MAXCUT problem, due to the existence of a large number of related decision variables and a large number of local optima, the optimization result of LNS is worse than Gurobi. Because ACP uses a two-step variable selection method to improve the correlation between variables to be optimized, and adaptively updates the block number to help jump out of the local optimum, the improvement is obvious compared with LNS on MAXCUT.

\subsection{Performances on Small-scale IPs}
\begin{table}[h]
\tiny
\begin{tabular}{|c|l|l|l|l|}
\hline
\textbf{Problem}                                                             & \textbf{\begin{tabular}[c]{@{}c@{}}SCIP \\ based\end{tabular}} & \textbf{Objective value}      & \textbf{\begin{tabular}[c]{@{}c@{}}Gurobi\\ based\end{tabular}} & \textbf{Objective value}      \\ \hline
\multirow{4}{*}{\begin{tabular}[c]{@{}c@{}}IS\\ (Maximize)\end{tabular}}     & SCIP                                                           & 1882.61 $\pm$ 23.74           & Gurobi                                                          & 2184.64 $\pm$ 14.21           \\
                                                                             & LNS                                                            & 2269.42 $\pm$ 12.02           & LNS                                                             & 2259.93 $\pm$ 7.67            \\
                                                                             & ACP                                                            & \textbf{2297.58 $\pm$ 13.53}  & ACP                                                             & 2291.69 $\pm$ 11.64           \\
                                                                             & ACP2                                                           & 2224.61 $\pm$ 166.78          & ACP2                                                            & \textbf{2293.66 $\pm$ 10.82}  \\ \hline
\multirow{4}{*}{\begin{tabular}[c]{@{}c@{}}MVC\\ (Minimize)\end{tabular}}    & SCIP                                                           & 3120.28 $\pm$ 38.69           & Gurobi                                                          & 2791.24 $\pm$ 20.13           \\
                                                                             & LNS                                                            & 2821.83 $\pm$ 34.87           & LNS                                                             & 2718.5 $\pm$ 15.9             \\
                                                                             & ACP                                                            & \textbf{2711.23 $\pm$ 28.55}  & ACP                                                             & 2686.44 $\pm$ 11.89           \\
                                                                             & ACP2                                                           & 2792.07 $\pm$ 150.51          & ACP2                                                            & \textbf{2682.87 $\pm$ 12.5}   \\ \hline
\multirow{4}{*}{\begin{tabular}[c]{@{}c@{}}MAXCUT\\ (Maximize)\end{tabular}} & SCIP                                                           & \textbf{10827.41 $\pm$ 28.34} & Gurobi                                                          & 10686.93 $\pm$ 53.98          \\
                                                                             & LNS                                                            & 10453.03 $\pm$ 38.3           & LNS                                                             & 10616.6 $\pm$ 57.51           \\
                                                                             & ACP                                                            & 10503.11 $\pm$ 52.69          & ACP                                                             & \textbf{11103.77 $\pm$ 99.91} \\
                                                                             & ACP2                                                           & \textbf{10827.41 $\pm$ 28.34} & ACP2                                                            & 10995.1 $\pm$ 61.76           \\ \hline
\multirow{4}{*}{\begin{tabular}[c]{@{}c@{}}SC\\ (Minimize)\end{tabular}}     & SCIP                                                           & 2507.64 $\pm$ 34.34           & Gurobi                                                          & 1800.23 $\pm$ 14.2            \\
                                                                             & LNS                                                            & 1815.14 $\pm$ 26.81           & LNS                                                             & 1693.23 $\pm$ 8.7             \\
                                                                             & ACP                                                            & \textbf{1656.86 $\pm$ 19.99}  & ACP                                                             & \textbf{1599.98 $\pm$ 8.85}   \\
                                                                             & ACP2                                                           & 1665.73 $\pm$ 23.19           & ACP2                                                            & 1604.34 $\pm$ 5.98            \\ \hline
\multirow{4}{*}{\begin{tabular}[c]{@{}c@{}}Read-world IP\\ (Maximize)\end{tabular}}    & SCIP                                                           & 18347.46 $\pm$ 22470.99           & Gurobi                                                          & 47378.52 $\pm$ 441.98          \\
                                                                             & LNS                                                            & 45694.59 $\pm$ 1098.18           & LNS                                                             & 46375.27 $\pm$ 281.97             \\
                                                                             & ACP                                                            & \textbf{46163.46 $\pm$ 1205.19}  & ACP                                                             & 47352.8 $\pm$ 427.65   \\
                                                                             & ACP2                                                           & 45747.41 $\pm$ 1096.76           & ACP2                                                            & \textbf{47379.18 $\pm$ 441.2}             \\ \hline
\end{tabular}
\caption{Comparison results for ACP and LNS with different subroutine solvers on different benchmark IPs with small-scale.}
\end{table}

\subsection{Performances on Medium-scale IPs}
\begin{table}[h]
\tiny
\begin{tabular}{|c|l|l|l|l|}
\hline
\textbf{Problem}                                                             & \textbf{\begin{tabular}[c]{@{}c@{}}SCIP \\ based\end{tabular}} & \textbf{Objective value}      & \textbf{\begin{tabular}[c]{@{}c@{}}Gurobi\\ based\end{tabular}} & \textbf{Objective value}      \\ \hline
\multirow{4}{*}{\begin{tabular}[c]{@{}c@{}}IS\\ (Maximize)\end{tabular}}     & SCIP                                                           & 18558.83 $\pm$ 53.08          & Gurobi                                                          & 21633.13 $\pm$ 86.11            \\
                                                                             & LNS                                                            & 19715.85 $\pm$ 32.75          & LNS                                                             & 21546.55 $\pm$ 345.81            \\
                                                                             & ACP                                                            & 21473.28 $\pm$ 272.57  & ACP                                                             & \textbf{22450.62 $\pm$ 383.8}          \\
                                                                             & ACP2                                                           & \textbf{21838.83 $\pm$ 38.04}          & ACP2                                                            & 22282.99 $\pm$ 629.87  \\ \hline
\multirow{4}{*}{\begin{tabular}[c]{@{}c@{}}MVC\\ (Minimize)\end{tabular}}    & SCIP                                                           & 31341.52 $\pm$ 80.79          & Gurobi                                                          & 28198.25 $\pm$ 32.95           \\
                                                                             & LNS                                                            & 30046.58 $\pm$ 146.88           & LNS                                                             & 28054.93 $\pm$ 32.82            \\
                                                                             & ACP                                                            & 28267.97 $\pm$ 254.05  & ACP                                                             & \textbf{26975.96 $\pm$ 24.43}          \\
                                                                             & ACP2                                                           & \textbf{27926.29 $\pm$ 140.61}          & ACP2                                                            & 26981.82 $\pm$ 32.18   \\ \hline
\multirow{4}{*}{\begin{tabular}[c]{@{}c@{}}MAXCUT\\ (Maximize)\end{tabular}} & SCIP                                                           & 0.99 $\pm$ 0.49               & Gurobi                                                          & 100211.7 $\pm$ 3103.54           \\
                                                                             & LNS                                                            & 94106.85 $\pm$ 615.97         & LNS                                                             & 95967.95 $\pm$ 1469.23           \\
                                                                             & ACP                                                            & \textbf{105299.09 $\pm$ 1481.79}          & ACP                                                             & 107857.56 $\pm$ 1501.38 \\
                                                                             & ACP2                                                           & 100770.17 $\pm$ 2092.36 & ACP2                                                            & \textbf{108244.69 $\pm$ 1374.97}           \\ \hline
\multirow{4}{*}{\begin{tabular}[c]{@{}c@{}}SC\\ (Minimize)\end{tabular}}     & SCIP                                                           & 25166.46 $\pm$ 37.69           & Gurobi                                                          & 17983.52 $\pm$ 40.89           \\
                                                                             & LNS                                                            & 24337.91 $\pm$ 265.94           & LNS                                                             & 17910.68 $\pm$ 67.92             \\
                                                                             & ACP                                                            & \textbf{20950.17 $\pm$ 526.61}  & ACP                                                             & 17243.17 $\pm$ 1290.13   \\
                                                                             & ACP2                                                           & 22445.23 $\pm$ 1024.17           & ACP2                                                            & \textbf{16388.32 $\pm$ 160.84}             \\ \hline
\multirow{4}{*}{\begin{tabular}[c]{@{}c@{}}Read-world IP\\ (Maximize)\end{tabular}} & SCIP                                                           & 0.0 $\pm$ 0.0                     & Gurobi                                                          & 454388.77 $\pm$ 2536.46          \\
                                                                                    & LNS                                                            & 454947.14 $\pm$ 2676.46           & LNS                                                             & 465150.57 $\pm$ 2475.65          \\
                                                                                    & ACP                                                            & \textbf{455054.64 $\pm$ 2746.93}  & ACP                                                             & \textbf{477051.41 $\pm$ 2613.14} \\
                                                                                    & ACP2                                                           & 454856.64 $\pm$ 2702.49          & ACP2                                                            & 472472.17 $\pm$ 2030.22          \\ \hline
\end{tabular}
\caption{Comparison results for ACP and LNS with different subroutine solvers on different benchmark IPs with medium-scale.}
\end{table}

\subsection{Performances on Large-scale IPs}
\begin{table}[h]
\tiny
\begin{tabular}{|c|l|l|l|l|}
\hline
\textbf{Problem}                                                                    & \textbf{\begin{tabular}[c]{@{}c@{}}SCIP \\ based\end{tabular}} & \textbf{Objective value}          & \textbf{\begin{tabular}[c]{@{}c@{}}Gurobi\\ based\end{tabular}} & \textbf{Objective value}         \\ \hline
\multirow{4}{*}{\begin{tabular}[c]{@{}c@{}}IS\\ (Maximize)\end{tabular}}            & SCIP                                                           & 7866.55 $\pm$ 55.47               & Gurobi                                                          & 215365.65 $\pm$ 113.5            \\
                                                                                    & LNS                                                            & 194733.92 $\pm$ 79.17             & LNS                                                             & 220368.34 $\pm$ 557.89           \\
                                                                                    & ACP                                                            & \textbf{207901.54 $\pm$ 1531.18}  & ACP                                                             & 227559.62 $\pm$ 57.26            \\
                                                                                    & ACP2                                                           & 196742.1 $\pm$ 116.6              & ACP2                                                            & \textbf{227636.51 $\pm$ 432.18}  \\ \hline
\multirow{4}{*}{\begin{tabular}[c]{@{}c@{}}MVC\\ (Minimize)\end{tabular}}           & SCIP                                                           & 490857.44 $\pm$ 164.24            & Gurobi                                                          & 283170.59 $\pm$ 317.06           \\
                                                                                    & LNS                                                            & 304860.39 $\pm$ 313.2             & LNS                                                             & 277870.09 $\pm$ 300.8            \\
                                                                                    & ACP                                                            & \textbf{291226.0 $\pm$ 1189.94}   & ACP                                                             & 271163.92 $\pm$ 337.0            \\
                                                                                    & ACP2                                                           & 301836.42 $\pm$ 70.59             & ACP2                                                            & \textbf{271087.22 $\pm$ 295.85}  \\ \hline
\multirow{4}{*}{\begin{tabular}[c]{@{}c@{}}MAXCUT\\ (Maximize)\end{tabular}}        & SCIP                                                           & 9.02 $\pm$ 1.24                   & Gurobi                                                          & 971138.1 $\pm$ 308.18            \\
                                                                                    & LNS                                                            & 553840.66 $\pm$ 45172.75          & LNS                                                             & 910662.36 $\pm$ 516.31           \\
                                                                                    & ACP                                                            & \textbf{829597.17 $\pm$ 46331.3}  & ACP                                                             & 1050400.97 $\pm$ 106.3           \\
                                                                                    & ACP2                                                           & 747434.19 $\pm$ 48697.78          & ACP2                                                            & \textbf{1053583.50 $\pm$ 473.43} \\ \hline
\multirow{4}{*}{\begin{tabular}[c]{@{}c@{}}SC\\ (Minimize)\end{tabular}}            & SCIP                                                           & 919264.06 $\pm$ 607.75            & Gurobi                                                          & 320034.15 $\pm$ 208.69           \\
                                                                                    & LNS                                                            & 584210.54 $\pm$ 158940.18         & LNS                                                             & 224435.73 $\pm$ 1253.96          \\
                                                                                    & ACP                                                            & \textbf{392082.53 $\pm$ 10206.15} & ACP                                                             & 208062.73 $\pm$ 5582.97          \\
                                                                                    & ACP2                                                           & 432368.51 $\pm$ 615.53            & ACP2                                                            & \textbf{201034.73 $\pm$ 196.14}  \\ \hline
\multirow{4}{*}{\begin{tabular}[c]{@{}c@{}}Read-world IP\\ (Maximize)\end{tabular}} & SCIP                                                           & 0.0 $\pm$ 0.0                     & Gurobi                                                          & 903119.36 $\pm$ 1588.62          \\
                                                                                    & LNS                                                            & 904387.43 $\pm$ 1560.34           & LNS                                                             & 924886.13 $\pm$ 1670.82          \\
                                                                                    & ACP                                                            & \textbf{909113.41 $\pm$ 2030.01}  & ACP                                                             & \textbf{948461.41 $\pm$ 4385.63} \\
                                                                                    & ACP2                                                           & 907155.48 $\pm$ 2199.08           & ACP2                                                            & 942235.76 $\pm$ 1861.36          \\ \hline
\end{tabular}
\caption{Comparison results for ACP and LNS with different subroutine solvers on different benchmark IPs with large-scale.}
\end{table}
\end{document}